\documentclass[12pt]{article}
\usepackage{amsmath}
\usepackage{xcolor}
\usepackage{authblk}
\usepackage{url, hyperref}
\usepackage{seqsplit}
\begin{document}

\newtheorem{definition}{Definition}[section]

\title{Finding a Widely Digitally Delicate Prime}
\author{Jon Grantham}
\affil{Center for Computing Sciences \\
Institute for Defense Analyses \\
Bowie, MD} 
\maketitle

\section{Introduction}

In a recent paper, Filaseta and Southwick \cite{FS} prove the remarkable theorem that a positive proportion of primes become composite when any digit is changed, {\it including leading zeros}. Such primes are known as {\bf widely digitally delicate primes}. This fact is all the more surprising when one considers that no such prime is explicitly known. Using an existing table of digitally delicate primes at the Online Encyclopedia of Integer Sequences \cite{oeis}, they show that there is no widely digitally delicate prime below $10^9$.

It is the aim of the current paper to remove that element of surprise. We do so by giving the first explicit example of a widely digitally delicate prime.

For a comprehensive history of related problems, see \cite{FS}. We will define a few terms for clarity. 

\subsection{Terminology}

\begin{definition}A {\bf covering system} is a collection of residue classes such that every integer is in at least one congruence class in the collection.
\end{definition}

\begin{definition}A {\bf digitally delicate prime} is a prime such that changing any one of its digits gives a composite. The digits under consideration do not include the leading zeros.
\end{definition}

\begin{definition}A {\bf widely digitally delicate prime} is a digitally delicate prime that also becomes composite when any one of its leading zeros is changed.
\end{definition}

\section{Finding the Covering System}

Filaseta and Southwick use a covering system to construct a congruence class such that any digitally delicate prime in that class is a widely digitally delicate prime. They then show that a positive proportion of primes in this congruence class are digitally delicate. 

Unfortunately, they do not explicitly construct the congruence class. They rely on the fact that ten specific unfactored composite numbers have at least two prime factors. In particular, any prime factor of the cyclotomic polynomial $\Phi_n(10)$ is a prime where $10$ has order $n$.  Any composite factor of $\Phi_n(10)$ (that is not a power) is divisible by at least two distinct primes $p_1,p_2$ where $10$ has order $n$. The unknown primes form two pieces of the covering system.

They note that it is possible to construct an explicit congruence class, but that the resulting numbers to be tested for primality would be around $20,000$ digits.

By modifying their original covering system, we only need to test numbers of slightly more than $4,000$ digits. We do not reproduce the congruence in
\cite{FS}, but recommend a careful reading of that paper to those who want to see the exact nature of the modifications described below.

The first step in their technique constructs a congruence class $a\bmod M$ such that for each digit $d=1,\ldots,9$, the quantity $a+d\times 10^k$ is guaranteed to be composite, for all $k$.
The recipe given for the digits $1$, $2$, $4$, $5$, $6$, $7$, $8$, and $9$ is very easy to follow --- all primes are given explicitly. The modifications required are to make sure that adding $3\times 10^k$ produces a composite.

The first observation is that the numbers of the form $\Phi_n(10)$ are factors of repunits. The practice of factoring numbers is
very popular on the Internet, as are repunits. Kamada \cite{kamada} runs a website that collects many of these factors. 

We use these factorizations to replace
$p_{242,1}$, $p_{242,2}$, $p_{275,1}$, $p_{275,2}$, $p_{363,1}$, $p_{363,2}$, $p_{396,1}$, $p_{396,2}$, $p_{484,1}$, $p_{484,2}$, $p_{605,1}$, $p_{605,2}$, $p_{726,1}$, $p_{792,1}$, $p_{1188,1}$, $p_{2420,1}$, $p_{4356,1}$, $p_{5808,1}$, and $p_{5808,2}$ with explicit prime values.

Kamada's website does not give us values for $p_{1210,1}$, $p_{1210,2}$, $p_{2904,1}$, and $p_{2904,2}$, so we have to improvise. Instead of $p_{1210,1}$, we use a third new factor of $\Phi_{605}(10)$, which we call $p_{605,3}$. Because $605$ is a divisor of $1210$, a congruence for $k$ modulo $605$ also holds modulo $1210$. Instead of $p_{1210,2}$, we use the remaining composite factor of $\Phi_{605}(10)$ after dividing the composite $C_{605}$ from \cite{FS} by the three new primes. Although knowing a prime factor would give a shorter overall congruence, a composite factor is equally valid.

Instead of $p_{2904,1}$, we use $p_{363,3}$, and instead of $p_{2904,2}$, we use $p_{726,3}$.

Filaseta and Juillerat \cite{FJ} have produced an alternate covering system, but the congruence class would be too large for the techniques in this paper.

\section{Implementation}

All of the code used in this computation is at \href{https://github.com/31and8191/delicate}{github.com/31and8191/delicate}.

In order to construct the congruence class for the prime, we put all of the congruence classes from \cite{FS} into a PARI/GP \cite{pari} script and then use the Chinese Remainder Theorem. This part of the computation takes a fraction of a second. We add the restriction $p\equiv 1\bmod 2$ so that we only search over odd numbers.

The next step is a C program using the PARI library. It proceeds through the congruence class, testing each number with PARI's {\tt ispseudoprime()} function. Once a number passes, the program tests each variation with a changed digit by using the same function. It is possible to inadvertently reject a digitally delicate number if a variation is a pseudoprime, but that is a small risk. Moreover, we are not concerned with false negatives.

That code ran for about $8$ hours on $50$ Xeon E5-2699 processors, each running at $2.3$ GHz, before producing a positive answer. Note that this answer is a widely digitally delicate {\bf probable} prime.

That brings us to the last step. I used the PARI/GP {\tt isprime()} function to prove the primality of the number. That proof took $30$ CPU-hours on Xeon processors described above. The prime contains $4032$ digits; it is given its own section below.

\section{Verification}

In order to verify that the prime is widely digitally delicate, I used PARI/GP to verify that each of the prime factors of the congruence modulus is a factor of $10^{439084800}-1$. Then I verified that $p+d\times 10^j$ for $1\le j \le 439084800$ and $1\le d \le 9$ has a non-zero GCD with the congruence modulus.

A previous version of this paper contained an incomplete verification, which did not catch an error in the construction of the prime. As a result, the prime was not verified to be
widely digitally delicate. The prime has been replaced. I am extremely grateful to the
referee for catching this error.

An ECPP primality certificate is at the GitHub repository.

\section{Other Computation}

I used the exact congruences provided in \cite{South} and similar code to provide examples of widely digitally delicate primes for bases $4$, $5$, $6$, $9$ and $11$.  Southwick had given examples for bases $2$ and $3$. Bases $7$ and $8$ should be equally easy.

I computed a list of digitally delicate primes up to $10^{11}$; none of them are widely digitally delicate.

\section{The Prime}

\seqsplit{285894570491987001178153724374587938515501125352188765520886436334395325908162183231168020433985595885849898174484619772705429763745991194664461100163727123429079686305371595295011006433565561943333249551496146898786776550562050563528909231406272849064430203150357126420677812739927895202546618565727359110480207958673019425654563382405130810590043829832715380016952742364731312668598733740964023055663765267200830877802707668653471933777180767950036688773110101833483505861901417331345780047986628791794326507501874968194285942890205589193464254902430887558565879364674586977542869103259950737623441567819481362009791661429525863338817583418337750636854201374491211035260749630474538478350982057067390182173675337406552432949290894282153403327225579837893941254410712722039794085468380534381131239387917463610086595526879843146781456368648816955679603677489665210301585644557371116185244779145233755095968806972088867036525332000563482661120311917367253346938362615371782983097701381566032284628055925636748898126820650464827690220782884190798343116909690410841541683928928146059651185336459496576809728038333262955417851244474541192580656810039197263444963081193148783858884674684291290222330733896006088695759175030396064034822228704608841906934859287094672381558955431945725037617470641744551971184997856292928496591534921211236168099241705439755139254869750171630873039529038801552200374073595158674067785012041849784559706868066251558959007352028060909364833937786587386985075068208343448492898929874743163115366206843020913389386172497510663145775351914248509567074373964353926550855085313243203493891126882173079296972161009951526577264968933125572307190814735375219620927527926344910610632478049454126284831311027582454324305015354771565731123702418064737323403806151773310815070975731441875702723456262422122995757282908371118151810622179472346353894805386163739494513263542196086802617326314422658695401574953618504325207626570369827335044357601565464406200866061462414874388427226428378855463876957838538058667545547011005526175955126923212015500245199711223649048445695487924834121120688806897224963006252750771969496500665634692331145254779980486343967092211316107065827196341907200951131788071359401349146624378022752052545844144930311500503905988265128637482488955981286529793467809095581440352363832112873230199489841731268453354502588069324055161268443331216962521399803300022233553354419767585266713647261688468914028856910813765524543565884834951623937277749823243669586887491792446589934309152360145402223702059870976018312422587314837330651273101910480339409326664560463527738695382754181288243632972126722859624445527255292587157441921530826113552602921728700347440467945156253246324775061989113520577721694133490784399906305246134471659941228999430202113778814079195332045833820656458901527522411086824910166211024837473058912892928937245474268883280684242480487874698030206629697898794434472627736433485220385368119684099619061236512191519875021429108929407668260565433218275587502476105227425889862121206865754501735824432651679406459045203787776692980263412366742788331072209322658634582697941520974990757805582353440188334311782884631555114175534159974173334421080734338168929817480812898846370530054327828885261092395358871603433955916853259029118665402592591212769522753087657181383889092152274888705253938497543051178220510833709854540377531040026335538883439006407548027704742980284105860308607796896327127214198076159411035112242982961440264096944328564609012135336998306466140086042893352680052097192600330626485072400319291412252578437918736963043218871598100051612110090345840283216887928981948760162910616377942967779894553056956511250045645042017691094927754292707440718212298955678804377942875527266376928166048118383050605417536635038030584921442558997485565729401000730292149168921659416381049130296671584907248290245314396708303820511731834810811383331428645024458864506105438802952519781329205899113620752868037776637274722119629626665263194401602076490693769774059953268451781114501333003}
\bibliographystyle{plain}
\bibliography{wide}

\end{document}